\documentclass[12pt,a4paper,doublesp]{amsart}
\usepackage[english]{babel}

\advance\topmargin -1.5cm
\advance\oddsidemargin -1.5cm
\advance\evensidemargin -1.5cm
\newcommand{\vect}{\operatorname{span}}

\newcommand{\be}[1]{\begin{equation} \label{eq:#1}}
\newcommand{\beq}{\begin{equation}}
\newcommand{\ee}{\end{equation}}

\newcommand{\bea}[1]{\begin{eqnarray} \label{eq:#1}}
\newcommand{\beaq}{\begin{eqnarray}}
\newcommand{\eea}{\end{eqnarray}}
\newcommand{\bean}{\begin{eqnarray*}}
\newcommand{\eean}{\end{eqnarray*}}

\newdimen\AAdi%
\newbox\AAbo%
\def\AArm{\fam0 \rm}%
\def\AAk#1#2{\setbox\AAbo=\hbox{#2}\AAdi=\wd\AAbo\kern#1\AAdi{}}%
\def\AAr#1#2#3{\setbox\AAbo=\hbox{#2}\AAdi=\ht\AAbo\raise#1\AAdi\hbox{#3}}%

\def\BBe{{\AArm I\!E}}%
\def\BBn{{\AArm I\!N}}%
\def\BBr{{\AArm I\!R}}%
\makeatletter

\@addtoreset{equation}{section}
\makeatother
\begin{document}

%%%%%%%%Declaracion de Teoremas %%%%%%%%%%

\newtheorem{thm}{Theorem}[section]

\newtheorem{lem}[thm]{Lemma}

\newtheorem{defn}[thm]{Definition}

\newtheorem{pro}[thm]{Proposition}

\newtheorem{cor}[thm]{Corollary}
\newtheorem{rem}[thm]{Remark}

%%%%%%%% FIN de Declaracion de Teoremas %%%%%%%%%%

%%%%%%%%     Titulo, autores y correos electronicos %%%%%%%%%%

\title[Iterative Methods for Linear Inverse Problems]{A Statistical view of  Iterative Methods for\\ Linear Inverse Problems}
\author[]{Ana K. Ferm\'{\i}n \and Carenne Lude\~na  }
\address{Escuela de Matem\'aticas, Facultad de Ciencias, UCV, Av. Los Ilustres, Los Chaguaramos, C\'odigo Postal 1020-A, Caracas -
 Venezuela. Telf.: (58)212-6051481.  \vskip .1in}
\address{Departamento de Matem\'aticas, IVIC, Carretera Panamericana KM. 11, Aptdo. 21827,
C\'odigo Postal 1020-A, Caracas -
 Venezuela. Telf.: (58)212-5041412. \vskip .1in}
\email{afermin@euler.ciens.ucv.ve} \email{cludena@ivic.ve}

%%% ----------------------------------------------------------------------

%%%%%%%%%%%%%%%%%  Agradecimientos   %%%%%%%%%%%%

\thanks{The authors would like to thank Projects \lq\lq Agenda Petr\'oleo" and \lq\lq ECOS-Nord" for their
financial support. }

%%% ---------------------------------------------------------------------

\begin{abstract}
In this article we study the problem of recovering the unknown
solution of a linear ill-posed problem, via iterative
regularization methods. We review the problem of
projection-regularization from a statistical point of view. A
basic purpose of the paper is the consideration of adaptive model
selection for determining regularization parameters. This article
introduces a new regularized estimator  which has the best
possible adaptive properties for a wide range of linear
functionals. We derive non asymptotic upper bounds for the mean
square error of the estimator and give the optimal convergence
rates.
\end{abstract}

%%% ----------------------------------------------------------------------


%%% ----------------------------------------------------------------------

\maketitle

%%% ----------------------------------------------------------------------

\section{Introduction}
\vskip .2in The area of mathematical inverse problems is quite
broad and involves the qualitative and quantitative analysis of a
wide variety of physical models. Moreover, a considerable number
of problems arising in different scientific and technical fields
belong to a class of ill-posed problems. For example,
geophysicists scan the earth's subsurface by recording arrival
times of waves reflected off different layers underneath the
surface, and try to determine a meaningful solution and to
understand which features in the solution are statistically
significant. \vskip .1in

From a statistical point of view, the problem can be seen as
recovering an unobservable signal $\tilde{f}$ based on
observations

\be{1} y(x_i)=A\tilde{f}(x_i)+\varepsilon_i, \ee where $A:F\to Y$
is some known compact linear operator defined over a  separable
Hilbert space $F$, with values in a separable Hilbert space $Y$
and $x_i,i=1,\ldots,n$ is a fixed observation scheme. We assume
that the observations $y(x_i) \in \BBr$ and that the observation
noise $\varepsilon_i $ are i.i.d. realizations of a certain random
variable $\varepsilon$. Throughout the paper, we shall denote
${\bf y}=(y(x_i))_{i=1}^n$. In this article we study the problem
of estimating $\tilde{f}$ using fixed iterative methods.\vskip
.1in

The best possible accuracy, regardless of any discretization and
noise corruption is determined by some \emph{a priori smoothness
assumption} on the exact solution $\tilde{f}$. Here, smoothness is
given in terms of some \emph{index function} $\eta$ on the
spectrum de $A^*A$ by
$$A_{\eta, \rho}=\{ f \in F, f= \eta(A^*A)\omega, \|\omega\| \leq
\rho\}$$ where $A_{\eta,\rho}$ is called a \emph{source
condition}. For classical Hilbert scales, the smoothness is
measured in terms of powers $\eta(t):= t^{\mu}$ with
$0\leq\mu\leq\mu_0$, $\mu_0 > 0$.\vskip .1in

In a deterministic framework, the statistical model \eqref{eq:1}
is formulated as the problem of finding the best-approximate
solution of $$Af=y,$$ in the situation where only perturbed data
$y^{\delta}$ are available with
$$\|y-y^\delta\|\leq \delta.$$ Here, $\delta$ is
called the noise level. It is important to remark that whereas in
this case consistency of the estimators depends on the
approximation parameter $\delta$, in (\ref{eq:1}) it depends on
the number of observations $n$. \vskip .1in

In general, the best $L^2$ approximation $A^{\dag}y$, where
$A^{\dag}$ is the Moore-Penrose (generalized) inverse of A, does
not depend continuously on the left-hand side $y$. We define the
Moore-Penrose inverse in an operator-theoretic way by restricting
the domain and range of $A$ in such a way that the resulting
restricted operator is invertible; its inverse will then be
extended to its maximal domain $\mathcal{D}(A^{\dag})=
\mathcal{R}(A) + \mathcal{R}(A)^{\bot}$, with $\mathcal{R}(A)$ the
range of the operator $A$ and $\mathcal{R}(A)^{\bot}$ the
orthogonal complement of the range of $A$.\vskip .1in

The inverse problems that we study in this article  are called
ill-posed problems because the operator $A$ is compact and
consequently equation \eqref{eq:1} can not be inverted directly
since $A^{-1}$ is not a bounded operator. Ill-posed problems are
usually treated by applying some linear regularization procedure,
often based on a singular value decomposition; see Tikhonov and
Arsenin in \cite{Tikhonov}. An interesting early survey of the
statistical perspective on ill-posed problems is studied in great
detail by O'Sullivan in \cite{osulliv}.\vskip .1in

In practice however \eqref{eq:1} is hardly ever considered.
Instead, we project the problem onto a smaller dimensional space
$Y_m$ of $Y$. This yields a sequence of closed subspaces $Y_m$
indexed by $m \in M_n$, a collection of index sets. Clearly, an
important problem is thus how to choose subspace $Y_m$ based on
the data. This can be done by selection of a cutoff point or by
threshold methods. Choosing the right subspace will be called
model selection.\vskip .1in

Sometimes this projection provides enough regularization to
produce a good approximate solution, but often additional
regularization is needed. Regularization methods replace an
ill-posed problem by a family of well-posed problems, their
solution, called regularized solutions, are used as approximations
to the desired solution of the inverse problem. These methods
always involve some parameter measuring the closeness of the
regularized and the original (unregularized) inverse problem,
rules (and algorithms) for the choice of these regularization
parameters as well as convergence properties of the regularized
solutions are central points in the theory of these methods, since
they allow to find the right balance between stability and
accuracy. The general principles of regularization for ill-posed
problems are known. In particular, such principles have been
established by A.N. Tikhonov. The literature on various
regularization methods based on these general principles is
extensive ( Engl, Hanke, and Neubauer \cite{Engl}, Gilyazov  and
Gol'dman \cite{Gilyazov} ).\vskip .1in

%In theoretical statistics, regularization is usually interpreted
%either in Bayesian terms. %Recent developments in the statistical
%literature offer promising new approaches to the fine-tuning of
%regularization techniques, particularly in the selection of
%regularization parameters.
%In applied statistics, regularization, often identified as \lq\lq
%penalty-based methods", \lq\lq soft thresholding", is associated
%primarily with nonparametric regression and density estimation,
%often referred to rather imprecisely as \lq\lq smoothing".

In statistic, regularization, is associated to penalty based
methods or thresholding methods or more generality to \lq\lq
smoothing" techniques. In applications, regularization offers a
unifying perspective for many diverse ill-posed inverse problems,
a wide range of problems concerned with recovering information
from indirect and usually noisy measurements, arising in
geophysics, tomography and econometrics. One of the most
important, but still insufficiently developed, topics in the
theory of ill posed problems is connected with iteration
regularization \cite{Gilyazov}; i.e, with the utilization of
iteration methods of any form for the stable approximate solution
of ill-posed problems. Iterative regularization  methods  tend to
be more attractive in terms of numerical cost and implementation,
but a number of open questions remain in their theoretical
analysis.\vskip .1in

In this article we propose an iterative regularized estimator for
linear \,ill-posed\, problems. Necessary conditions for
convergence are established. These conditions connect the choice
of the regularization parameter (i.e., the iteration index) with
the projection dimension. Moreover, we prove that the iterative
regularized estimator is optimal in the sense that the estimator
achieves the best rate of convergence among all the regularized
estimators. A recent work in this direction is developed  by
Loubes and Lude\~na in \cite{Jmi}, which discusses the problem of
estimating inverse nonlinear ill-posed problems with different
types of complexity penalties leading either to a model selection
estimator or to a regularized estimator. \vskip .1in

%At first glance, there can appear to be a lot of work associated
%with the selection of a good regularization parameter, and many
%algorithms proposed in the literature tend to be hard to implement
%are needlessly complicated. The choice of the regularization
%sequence is here crucial. Our investigation about properties of
%regularization, for ill-posed problems, using iterative methods,
%begins with the necessity of finding a good stopping criteria for
%terminating the iteration.

The choice of the regularization sequence is here crucial, and a
lot of work associated with the selection of a good regularization
parameter  can be found in the literature \cite{Engl-Grever},
\cite{Kilmer}. When using iterative methods the problem is finding
a good stopping criteria for terminating the iteration procedure.
In this article we will use tools developed in the context of
model selection via penalization, \cite{Birge},\cite{Barron},
based on the use of concentration inequalities. \vskip .1in

Our article is organized as follows. Section 2 presents basic
assumptions and a statement of the discretized inverse problem. In
section 3 we discuss regularization methods and we prove
consistency of the estimator when the regularization parameter is
known. In section 4 we present our main result, prove optimality
of an adaptive regularized estimator and give its rate of
convergence. Finally, in the last section we introduce
regularization by iterative methods for the solution of inverse
problems and provide some examples to explain the properties of
iterative regularization methods.
%%% ----------------------------------------------------------------------
\vskip .3in

\section{Preliminaries}\label{prelim}
\vskip .2in

\subsection{Formulation of the problem and basic assumptions } $ $ \vskip .2in

We assume  that the inverse problem is given by
$$y=A\tilde{f}+\varepsilon.$$
where $\varepsilon$ is a centered random variable satisfying the
moment condition $\BBe (|\varepsilon|^p/\sigma^p) \le p ! /2$  and
$\BBe (\varepsilon^2)=\sigma^2$. \vskip .1in

We also need some notations concerning the fixed design settings,
$x_i,i=1,\ldots,n$. Define the empirical measure:
$$P_n = \frac{1}{n} \sum_{i=1}^n \delta_{x_i}.$$
and the associated empirical norm
$$\| y \|_{n}^2=\| y \|_{P_n}^2 = \frac{1}{n} \sum_{i=1}^n (y(x_i))^2  $$
as well as the empirical scalar product $$<y,\varepsilon
>_n=\frac{1}{n} \sum_{i=1}^n \varepsilon_i y(t_i).$$
\vskip .1in We assume $A : F \to Y$, $F,\,Y$ separable Hilbert
spaces. Let $\langle \,,\, \rangle$ stands for the inner product
defined over $F$. We assume that the range of the operator $A$,
$\mathcal{R}(A),$ is closed in the sense that $\tilde{f} \in
\mathcal{N}(A)^{\bot}$, where $\mathcal{N}(A)^{\bot}$ is the
orthogonal complement of the null space of the operator $A$.\vskip
.1in
%The prove existence and uniqueness is an immediate consequence of the condition:
%the solution exists if and only if the operator equation
%Af=\Pi_{\bar{R(A)}}y is solvable.

With this notation let $J(f)= \|y - A f\|_n^2$ the quadratic risk
function.
%We generalize the inverse problem and instead seek the problem of
%risk minimization. The goal is to minimize the quadratic risk
%function $J(f)= \|y - A f\|_n^2$.
We will denote by  $\hat{f}$ the function that minimizes the risk
(which may not be unique), defined as

\be{2} \hat{f}= {\rm arg}\min_{f \in F} J(f),\ee where the minimum
is taken over all functions from $F$ to $Y$.\vskip .1in

The solution of the problem $\min_{f \in F} J(f)$ exists if and
only if $f$ is a solution of the normal equation
\be{ecnormal}A^*Af=A^*y\ee where $A^*: Y \to F$ is the adjoint
operator of $A$ (introduced via the requirement that for all $f
\in F$ and $y\in Y$, $\langle Ax, y \rangle_n = \frac{1}{n}
\langle x,A^*y \rangle$ holds).\vskip .1in

It is important to remark that the operator $A^*$ actually depends
on the observation sequence $x_i, \, i=1,\ldots n$. If $Y$ is
generated by  $\{\varphi_j\}_{j=1}^m$ and is such that this basis
is orthonormal with respect to the $L^2(P_n)$ norm over $Y$, and
$A$ is the identity then $A^* =\frac{1}{n}(
\varphi_j(x_i))_{i,j}$. \vskip .1in

It is necessary to mention that the convergence rates can thus be
given only over subsets of $F$, i.e., under  a-priori assumptions
on the exact solution $\tilde{f}$. We will formulate such a-priori
assumptions, encountered typically in the inverse problem
literature, in terms of the exact solutions by considering subsets
of $F$ given by some \emph{source condition} of the form

$$A_{\mu,\rho}=\{ f \in F, f=
(A^*A)^{\mu} \omega, \|\omega\| \leq \rho\}$$ where
$0\leq\mu\leq\mu_0$, $\mu_0 > 0$ and use the  notation

\be{3} A_{\mu}=\bigcup_{\rho>0}
A_{\mu,\rho}=\mathcal{R}((A^*A)^{\mu})\ee \vskip .1in

These sets are usually called source sets, $f \in A_{\mu, \rho}$
is said to have a source representation. The requirement that $f$
be in $A_{\mu, \rho}$ can be considered as an smoothness
condition. \vskip .3in
%%% ----------------------------------------------------------------------

\subsection{Projection methods}$ $ \vskip .2in

For numerical calculations, we have to approximate the space $F$
by a finite-dimensional subspace. Estimating over all $F$ is in
general not possible. One approach in this direction is
\emph{regularization by projection}, where the regularization is
achieved by a finite-dimensional approximation through projection.
\vskip .1in

Let $M_n$ be a collection of index sets ($m \in M_n, m =
\{j_1\ldots, j_{d_m}\}$). We give a sequence $Y_1\subset Y_2\ldots
\subset Y_m \ldots \subset Y$ whose of union is dense in $Y$. We
assume
$$dim(Y_m)=d_{m}.$$

Let $\Pi^n_{Y_m}$ be the orthogonal projector in the empirical
norm over the subspace $Y_m$ and let $A_m= \Pi^n_{Y_m}A$. Define
$F_m= A^* Y_m$,  with $A_m^* : Y_m \to F$ the adjoint operator of
$A_m$, and $\Pi_{F_m}$ to be the orthogonal projector onto the
subspace $F_m$. Then, by construction $$\Pi_{F_m}=(\Pi^n_{Y_m}A)^+
\Pi^n_{Y_m} A.$$

Thus, we shall assume that data are give through an orthogonal
design, corresponding to an orthogonal projection $\Pi^n_{Y_m}$ as

\be{new} \Pi^n_{Y_m} y = \Pi^n_{Y_m} A \tilde{f} + \Pi^n_{Y_m}
\varepsilon. \ee

%We shall assumed, that data are given through an orthogonal
%design, corresponding to an orthogonal projection. Let $\Pi_{F_m}$
%the orthogonal projector onto the subspace $F_m$ and $\Pi^n_{Y_m}$
%the orthogonal projector, in the empirical norm, over the subspace
%$Y_m$. \vskip .1in

With this notation we have that the best-approximate $L^2$
solution has the expression
$$\Pi_{F_m}\tilde{f} = A_m^{\dag}y_m.$$  for $y_m = \Pi_{Y_m}^n y $ in the domain of
$A_m^{\dag}$. In the following  we shall denote
$\tilde{f}_m=\Pi_{F_m}\tilde{f}$.\vskip .1in

%We have that for projection without regularization the choice of
%$F_m$ and of $A_m$ has many advantages. For instance, we have that
%$A_m^{\dag}y=\Pi_{F_m}\tilde{f}$, where $\tilde{f}=A^{\dag}y$.
%Thus, in the noise free case $A_m^{\dag}y$ is the best possible
%approximation over $F_m$ and there in no need for further
%regularization. However,

Our goal is to find the solution of the equation \eqref{eq:1} in
the finite-dimensional subspace $F_m$ of $F$. We have that for
projection without regularization the choice of $F_m$ and of $A_m$
has many advantages. For noisy data and severely ill-posed
problems the dimension of the subspace has to be rather low to
keep total error estimate small, since for these problems the
smallest singular value of $A_m$ decreases rapidly as $m$
increases. To be able to use larger dimensions we have to combine
the projection method with additional regularization methods,such
as iterative methods \cite{Engl},\cite{Gilyazov}. \vskip .1in

%If $Y_m$ is generated by some orthonormal basis $\{\varphi_j\}_{j
%\in m}$, $Y_m = \vect \{\varphi_j\}_{j \in m}$, with respect to
%the $L^2(P_n)$ norm over $Y$, then
%$$\Pi_{Y_m}^n y=\sum_{j=1}^{d_m} y_{j,n}\,\varphi_j ,$$ where
%$y_{j,n} =\,<\Pi_{Y_m}^n y,\varphi_j>_n$ are the solution to the
%projection problem under the empirical measure $P_n$. Set
%$\Sigma_m=[\varphi_j(t_i)]_{{(i)_1^n,}\,{j \in m}}$ to be the Gram
%matrix associated to basis $\left(\varphi_j \right)_{j \in m}$
%over $Y_m$. Thus, $$y_{j,m}= (\Sigma_m^t \Sigma_m)^{-1} \Sigma_m^t
%(y(t_1),\ldots, y(t_n)).$$
\vskip .3in
%%% ----------------------------------------------------------------------

\subsection{Singular value decomposition}$ $ \vskip .2in

%The primary difficulty with linear ill-posed problems is that the
%inverse image is undetermined due to small (or zero) singular
%values of $A_m$. Usually that is not significant for the large
%singular values, but may lead to ambiguity in the small singular
%values so that we do not know if they are small or zero. \vskip
%.1in

As often $A_m$ is not of full rank, the singular value
decomposition (SVD) of the operator $A_m$ is then a useful tool.
Let $(\sigma_j; \phi_j, \varphi_j)_{j \in m}$ be a singular system
for a linear operator $A_m$, that is, $A_m \phi_j = \sigma_j
\varphi_j$ and $A_m^* \varphi_j = \sigma_j \phi_j$; where
$\{\sigma_j^2\}_{j \in m}$ are the non–zero eigenvalues of the
self–adjoint operator $A_m^* A_m$ (and also of $A_m A_m^*$),
considered in decreasing order. Furthermore, $\{\phi_j\}_{j \in
m}$ and $\{\varphi_j\}_{j \in m}$ are a corresponding complete
orthonormal system of eigenvectors of $A_m^*A_m$ and $A_m A_m^*$,
respectively. For general linear operators with an SVD
decomposition, we can write

\be{5} A_m f = \sum_{j \in m} \sigma_{j} \langle f , \phi_j\rangle
\varphi_j  \ee

\be{6} A_m^* y_m = \sum_{j \in m} \sigma_{j} \langle y_m ,
\varphi_j\rangle \phi_j.\ee  \vskip .1in

For $y_m$  in the domain of  $A_m^{\dag}$,
$\mathcal{D}(A_m^{\dag})$, the best-approximate $L^2$ solution
hast the expression
$$A_m^{\dag}y_m = \sum_{j \in m} \frac{\langle y_m ,\varphi_j\rangle
}{\sigma_j} \phi_j = \sum_{j \in m} \frac{\langle A_m^*y_m
,\phi_j\rangle}{\sigma_j^2} \phi_j.$$

Note that for large $j$, the term $1/\sigma_j$ grows to infinity.
Thus, the \emph{high frequency errors} are strongly amplified. We
will asume that $\sigma_j = \mathcal{O}(j^{-p})$ for some $ p
> 1/2$, which is clearly related to the ill-posedness of the
operator $A$ and the approximation properties of $Y_m$. For the
construction and analysis of regularization methods, we will
require some general  notation for functions of the operators
$A_m^* A_m$  and  $A_m A_m^*$ . \vskip .1in

%Also, the faster the decay of the $\sigma_j = \mathcal{O}(j^{-p})$
%(with $ p
%> 1/2$). This rate depends on the ill-posedness of the operator and the
%approximation properties of $F_m$. In some cases these are known
%precisely and in others they can be deduced from the properties of
%the solution $\tilde{f}$. One such case is the following source
%assumption encountered typically in the inverse problems
%literature \cite{Engl}. Regularization methods now are techniques
%that can handle these problems. For the constructions and analysis
%of regularization methods, we will need la notion of spectral
%functions of selfadjoint operators. \vskip .1in

%We denote spectral functions \ of \ the \ operators $A_m^* A_m$ \
%and \ $A_m A_m^*$ \ by \ $E_{\lambda}$ and $H_{\lambda}$,
%respectively.\  The spectral family \ $E_{\lambda}$ is called \
%\emph{spectral \ decomposition} of \ $A_m^* A_m$.\vskip .1in

%For all $\lambda$,
%$$E_{\lambda}(\cdot) = \sum_{\sigma^2 < \lambda_{j}, j \in
%m}\langle \cdot,\phi_j \rangle \phi_j$$ is an orthogonal projector
%and projects onto
%$$F_{\lambda}= \vect \{\phi_j \,|\, j \in m,\, \sigma^2 < \lambda\}.$$\vskip .1in

Let $E_{\lambda}$ be the spectral decomposition of $A_m^* A_m$
given by

$$E_{\lambda}(\cdot) = \sum_{\sigma^2 < \lambda_{j}, j \in
m}\langle \cdot,\phi_j \rangle \phi_j$$ and $H_{\lambda}$ the
spectral decomposition of $A_m A_m^*$. Then $E_{\lambda}$ is an
orthogonal projector and projects onto $$ \vect \{\phi_j \,|\, j
\in m,\, \sigma^2 < \lambda\}.$$\vskip .1in

Since $(\sigma_j^2; \phi_j)$ is an eigensystem for the selfadjoint
compact operator $A_m^*A_m$,
$$ A_m^*A_m f = \sum_{j \in m} \sigma^2_{j} \langle f ,
\phi_j\rangle \phi_j $$ holds, which will be written (using the
definition of the integral below) as
$$ A_m^* A_m f = \int \lambda \,dE_{\lambda} f $$
for $f \in \mathcal{D}(A_m)$. Here the limits of integration could
be 0 and $\|A_m\|^2 + \epsilon$ for any $\epsilon > 0$. We
sometimes omit the limits of integration. \vskip .1in

This, motivates the definition \be{4} G(A_m^*A_m):= \int
G(\lambda) dE_{\lambda} := \sum_{\sigma^2 = \lambda_{j}, j \in m}
G(\sigma_j^2) \langle \cdot , \phi_j\rangle \phi_j \ee of a
(piecewise) continuous function $G$ of a selfadjoint linear
operator on $F_m$. If $A^*A$ is continuously invertible, then
$(A^*A)^{-1}=\int \frac{1}{\lambda}\,dE_{\lambda}$.\vskip .1in

In this case the best-approximate $L^2$ solution, for $y_m$ in the
domain of $A^{\dag}_m$, can be characterized by the equation

\be{7}\tilde{f}_m = A_m^{\dag} y_m = \int\, \frac{1}{\lambda}\,
dE_{\lambda}A_m^* y_m.\ee

If $G(A_m^* A_m)$ is defined via \eqref{eq:4}, then for $f\in
\mathcal{D}(G_1(A_m^*A_m))$ and $g \in \mathcal{D}(G_2(A_m^*A_m))$

\be{8}\langle G_1(A_m^* A_m) f, G_2(A_m^*A_m) g \rangle = \int
G_1(\lambda)G_2(\lambda)\,d \langle E_{\lambda}f , g \rangle\ee
and \be{9}\| G(A_m^* A_m) f \|^2 = \int G(\lambda)\, d \|
E_{\lambda} f \|^2.\ee

\vskip .1in

The source set, $A_{\mu}$ \eqref{eq:3}, can be characterized via
the singular values as follows:
$$ (A_m^*A_m)^{\mu} \omega = \int \lambda^{\mu} dE_{\lambda} \omega =
\sum_{j \in m} \sigma_j^{2 \mu} \langle \omega, \phi_j\rangle
\phi_j $$

%%% ----------------------------------------------------------------------
\vskip .3in

\section{Regularization methods }\label{RM}\vskip .2in

After the general considerations of the last section, we now
explain the construction of a regularization method for the
important special case of selfadjoint linear operators. The basic
idea for deriving a regularization method is to replace the
amplification factors $1/\lambda_j$ by a filtered version
$Q(\lambda_j,\alpha)$, where the filter function is a piecewise
continuous, nonnegative and nonincreasing function of $\lambda$ on
the segment $[0,\|A_m\|^2]$ for a regularization parameter $\alpha
> 0$. The assumptions over the regularizing
coefficients $Q_{\alpha}(\lambda)$ are technical and are given in
order to control fluctuations over set $[0,\|A_m\|^2]$.\vskip .1in

The filter family \ $\{Q(\lambda_j,\alpha)\}_{j\in m}$ \
approximates the function \ $\lambda^{-1}$ for  $\alpha  \to
\infty$. Intuitively, a regularization on \ $A_m^{\dag}$ \ should
then be the replacement of the ill conditioned operator \
$A_m^{\dag}$ \ by a \, family $\{R(\lambda_j,\alpha)\}_{j\in m}:
Y_m \to F_m$ of continuous operators. Throughout all the article,
we shall denote $\{Q(\lambda_j,\alpha)\}_{j\in m}$ and
$\{R(\lambda_j,\alpha)\}_{j\in m}$ by $Q_{\alpha}$ and
$R_{\alpha}$, respectively. Obviously, for all $\alpha > 0,$
$R_{\alpha}$ is bounded.  \vskip .1in

As the approximation of $\tilde{f_m}$, we then take

$$f_{m,\alpha}=Q_{\alpha}(A_m^*A_m) A_m^* y_m = R_{\alpha} y_m,$$ where
$R_{\alpha}:=\int Q_{\alpha}(\lambda) dE_{\lambda}A_m^*$.

\vskip .1in \vskip .1in
\begin{rem}
Note that with the above notation \be{88}f_{m,\alpha} =
R_{\alpha}y_m= R_{\alpha}\Pi_{Y_m}^n A \tilde{f} +
R_{\alpha}\Pi_{Y_m}^n \varepsilon.\ee Also that we can write
$$R_{\alpha}=Q_{\alpha}(A_m^*A_m)A_m^*.$$
\end{rem}
%\begin{rem}
%The assumptions over the regularizing coefficients
%$Q_{\alpha}(\lambda)$ are technical and are given in order to
%control fluctuations over set $[0,\|A_m\|^2]$.
%\end{rem}
\vskip .1in

%Thus, a regularization method consists of a regularization
%operator and a parameter choice rule which is convergent in the
%sense that, if the regularization parameter is chosen according to
%that rule, then the regularized solutions converge. In a
%deterministic framework, a requirement for regularization
%parameter $\alpha$ is that, if the noise level $\delta$ tends to
%0, the regularized solution $f_{m,\alpha}$ should then to
%$A^{\dag}y$. In this article, we would choose $\alpha$ from the
%data in such a way that optimal rates are maintained. This choice
%must also not depend on a priori regularity assumptions. The
%following question now arise: Under that condition ${R_{\alpha}}$
%it is an regularization operator for $A_m^{\dag}$?. \vskip .1in

%This question is discussed extensively in the regularization
%theory \cite{Engl-H}, \cite{Engl}.

The \,next \,theorem \, gives \, conditions \,under \,which \,the
first term in \eqref{eq:88} converges to $\tilde{f}_m =
\Pi_{F_m}\tilde{f}$. The proof follows that of \cite{Engl}, but we
include it for the sake of completeness.\vskip .2in

\begin{thm}\label{teo1}
Let, for all $\alpha >0$, $Q_{\alpha}:[0,\|A_m\|^2]\rightarrow
\BBr$ be a piecewise continuous and nonincreasing function of
$\lambda$ on the segment $[0,\|A_m\|^2]$. Assume also that there
is a $C > 0$ such that
$$|\lambda Q_{\alpha}(\lambda)| \leq C,$$
and
$$\lim_{\alpha \to \infty} Q_{\alpha}(\lambda) = \frac{1}{\lambda}$$
for all $\lambda \in (0,\|A_m\|^2).$ Then, for all $y \in
D(A_m^{\dag})$, \vskip .1in
$$\lim_{\alpha\rightarrow\infty} Q_{\alpha}(A_m^*A_m)A_m^* \Pi_{Y_m}^n A
\tilde{f} = \tilde{f}_m$$ holds with $\tilde{f}_m =
A_m^{\dag}y_m.$
\end{thm}\vskip .1in
\begin{rem} In order to assume convergence as $\alpha \to \infty$, it is necessary to choose $Q_{\alpha}$ such that it approximates
$1/\lambda$ for all $\lambda \in (0,\|A_m\|^2]$. Also, note that
the condition $|\lambda Q_{\alpha}(\lambda)| \leq C$ implies that
$\|A_m R_{\alpha}\|=\|A_m A_m^* Q_{\alpha}(A_m^*A_m) \| \leq C$,
i.e, $\|A_m R_{\alpha}\|$ is uniformly bounded.
\end{rem}

\vskip .1in

\begin{proof}
As in \cite{Engl}, if $\tilde{f}_m$ is defined by \eqref{eq:7},
then by \eqref{eq:ecnormal} the residual norm has the
representation
$$\|\tilde{f}_m-Q_{\alpha}(A_m^*A_m)A^*_m A_m \tilde{f}\|^2=
\|(I-Q_{\alpha}(A_m^*A_m)A_m^*A_m)\tilde{f_m}\|^2$$ \vskip .1in
From the formula \eqref{eq:9}, it follows that
$$\|\tilde{f}_m-Q_{\alpha}(A_m^*A_m)A^*_m A_m \tilde{f}\|^2 = \int_0^{\|A_m\|^2 +} (1-\lambda Q_{\alpha}(\lambda))^2\, d
\|E_{\lambda}\tilde{f}_m\|^2.$$\vskip .1in

Since $(1-\lambda Q_{\alpha}(\lambda))^2$  is bounded by the
constant $(1+C)^2$, which is integrable with respect to the
measure $d \|E_{\lambda}\tilde{f}_m\|^2,$  then by the Dominated
Convergence Theorem, \be{10}\lim_{\alpha \to \infty}
\int_0^{\|A_m\|^2 +} (1-\lambda Q_{\alpha}(\lambda))^2\,d
\|E_{\lambda}\tilde{f}_m\|^2 = \int_0^{\|A_m\|^2 +} \lim_{\alpha
\to \infty} (1-\lambda Q_{\alpha}(\lambda))^2\,d
\|E_{\lambda}\tilde{f}_m\|^2.\ee \vskip .1in Since for $\lambda
>0$, $\lim_{\alpha \to \infty}(1-\lambda Q_{\alpha}(\lambda))=0 $
then the integral on the right-hand side of \eqref{eq:10}  equals
to $0$. On the other hand, if $\lambda=0$, $\lim_{\alpha \to
\infty}(1-\lambda Q_{\alpha}(\lambda))=1$ then the equation
\eqref{eq:10} has the form \be{11}\lim_{\alpha \to \infty}
\int_0^{\|A_m\|^2 +} (1-\lambda Q_{\alpha}(\lambda))^2\,d
\|E_{\lambda}\tilde{f}_m\|^2 = \lim_{\lambda \to 0^+}
\|E_{\lambda}\tilde{f}_m\|^2 - \|E_{0}\tilde{f}_m\|^2\ee which is
equal the jump of $\|E_{\lambda}\tilde{f}_m\|^2$ at $\lambda = 0$.
Since $\tilde{f}_m \in \mathcal{N}(A_m)^{\bot}$,  the term on the
right-hand side of \eqref{eq:11} equals to 0. Thus, $R_{\alpha}
A_m \tilde{f}$ converges to $\tilde{f}_m$ as $\alpha \to \infty$
for $y_m \in D(A_m^{\dag}),$ which ends the proof.

\end{proof}
\vskip .2in

Let $Tr(B)$ the trace of the selfadjoint operator $B^tB$ for any
square matriz $B$, which is defined by
$$Tr(B)=\frac{1}{n}\sum_{j\in m} b_j$$ for $b_j$ eigenvalues of $B^tB$.

\vskip .1in

We then have the following result, \vskip .2in

\begin{thm}\label{teo2}
Let  $Q_{\alpha}$ be as in theorem \ref{teo1}. Let $\mu, \rho > 0$
and let $\omega_{\mu}: (0, \alpha_0) \to R$ be such that for all
$\alpha \in (0,\alpha_0)$ and $\lambda \in [0,\sigma_1^2]$

$$\sup_{0\leq \lambda \leq \sigma_1^2} \lambda^{\mu} |1 - \lambda Q_{\alpha}(\lambda)| \leq \omega_{\mu}(\alpha)$$
holds. Then for $\tilde{f} \in A_{\mu, \rho},$ the following
inequality holds true
\be{99}\BBe \|\tilde{f}_m - f_{\alpha,m}\|^2
\leq 2\,\omega_{\mu}(\alpha)^2\,\rho^2 +
2\,\sigma^2\,Tr(Q_{\alpha}^2(A_m^*A_m)A_m A_m^*).\ee
\end{thm}
\vskip .1in

\begin{proof}
The proof of this inequality is based on the definition of the
estimator $f_{m,\alpha}$ and on the assumptions over this
function. We have that the $L^2-$norm of the difference between
the regularized function and the true data function can be bounded
by

\be{a}\BBe\| \tilde{f}_m-f_{m,\alpha} \|^2 \leq 2\,
\BBe\|\tilde{f}_m-R_{\alpha}A_m\tilde{f }\|^2 + 2\, \BBe\|
R_{\alpha}A_m \tilde{f} - f_{m,\alpha}\|^2\ee where
$f_{m,\alpha}=R_{\alpha}y_m $. \vskip .1in

This is the typical bias-variance decomposition. The first term on
the right-hand side is an approximation error, which corresponds
to the bias, whereas the second term, variance, is a stability
bound on the regularizing operator $R_{\alpha}$. Note that by the
Theorem \ref{teo1}, the first term in \eqref{eq:a} goes to 0 if
$y_m \in \mathcal{D}(T_m^{\dag})$.\vskip .1in

Let $\omega \in F_m$ with $\|\omega\|\leq \rho$. Since
$\tilde{f}\in F_m$ then $\Pi_{F_m}\tilde{f}=(A_m^*A_m)^{\mu}\omega
$. On the other hand, $\lambda^{\mu}\sup_{ \lambda} |1 - \lambda
Q_{\alpha}(\lambda)| \leq \omega_{\mu}(\alpha)$,
 then the first term in this
equation can be bounded by

\begin{align*}
\BBe\| \tilde{f}_m-R_{\alpha}A_m\tilde{f}\|^2 & = \BBe\| \tilde{f}_m-Q_{\alpha}(A_m^*A_m)A_m^*A_m\tilde{f}_m\|^2\\
& = \BBe\| (I-Q_{\alpha}(A_m^*A_m)A_m^*A_m) \tilde{f}_m \|^2\\ &
\leq \omega_{\mu}^2 \, \rho^2.
\end{align*}

In order to control the term corresponding to the variance  we
used that the data perturbation is white noise. Thus,

\begin{align*}
    \BBe \|R_{\alpha}A_m\tilde{f}  - f_{m,\alpha}\|^2
    & = \BBe \langle  \varepsilon,(Q_{\alpha}(A_m^*A_m)A_m^*)^*Q_{\alpha}(A_m^*A_m)A_m^* \varepsilon
    \rangle\\
    & = \BBe \langle \varepsilon,Q_{\alpha}^2(A_m^*A_m)A_m A_m^*  \varepsilon
    \rangle\\
    & =  \sigma^2 \,Tr(Q_{\alpha}^2(A_m^*A_m)A_m A_m^*)
\end{align*}
which yields the desired result.
\end{proof}

\vskip .1in

The next result will be useful when studying iterative methods.
\vskip .2in

\begin{thm}\label{teo3}
Let  $Q_{\alpha}$ be as in theorem \ref{teo1}. Assume also that
$Q_{\alpha}$ is continuously differentiable and that the function
$|1-\lambda Q_{\alpha}(\lambda)|' |\lambda
Q_{\alpha}(\lambda)-1|^{-1}$ doest not decrease. Then the
estimates are valid
$$\sup_{0\leq \lambda \leq \sigma_1^2} |Q_{\alpha}(\lambda)| = Q_{\alpha}(0),$$
and
$$\sup_{0\leq \lambda \leq \sigma_1^2} \lambda^{\mu}|1-\lambda
Q_{\alpha}(\lambda)|< \mu^{\mu} (\mu +1)^{-1}\omega_{\mu}(\alpha)
$$ where $\omega_{\mu}(\alpha)= Q_{\alpha}(0)^{-\mu}.$

\end{thm}
\vskip .1in
\begin{proof}
The proof can be carried out by standard techniques. A proof of
this result can be found in \cite{Gilyazov}.
\end{proof}

%%% ----------------------------------------------------------------------
\vskip .3in

\section{Rates of convergence for the regularized estimator }\label{conv}
\vskip .2in

In any regularization method, the regularization parameter
$\alpha$ plays a crucial role. For choosing the parameter, there
are general methods of parameter selection. For example, the
Discrepancy Principe \cite{Morozov}, Cross-Validation \cite{cross}
and the L-curve \cite{Engl-Grever}. They differ in the amount of a
priori information required as well as in the decision criteria.
The appropriate choice of regularization parameter is a difficult
problem. We would like too choose $\alpha$, based on the data in
such a way that optimal rates are maintained. This choice should
not depend on a priori regularity assumptions.\vskip .1in

Our goal is to introduce  adaptive methods in the context of
statistical inverse problems. In this section we introduce our
adaptive estimator, for a fixed $m=m_0$. We choose $m_0$ such that
$\|\tilde{f}-\Pi_{F_{m_0}}\tilde{f}\|^2$ satisfies the optimal
rates with high probability since we know
$$\|\tilde{f}-\Pi_{F_{m_0}}\tilde{f}\|^2 < \|I-\Pi_{Y_{m_0}}\|^{4\mu} = \mathcal{O}(d_{m_0}^{-4\mu p}) $$
for a certain $p$ and $0<\mu\leq 1/2$. It is satisfied if the
dimension of the set is such that

\be{condition}d_{m_0}\geq n^{\frac{1}{2p +1}}.\ee \vskip .1in

%::::::::::::::::::::::::::::::::::::::::::::::::::::::::::::::::::::::::
%NOTA: NO HAY QUE PONER d_{m_0}\geq \alpha_n^{\frac{1}{4\mu p + 2p +1}}
%ya que dependeria de \mu y mu es desconocido solo sabemos que es mayor
%que 0 y menor que 1/2. Entonces como 2p+1 < 2p + 1 + \varepsilon
%entonces para \varepsilon = 2\mu p tenemos que se cumple que
%d_{m_0}\geq \alpha_n^{\frac{1}{2p + 1}} \geq n^{\frac{1}{4\mu p + 2p +1}}
%::::::::::::::::::::::::::::::::::::::::::::::::::::::::::::::::::::::::

This leads to the rate

$$\|\tilde{f}-f_{m,\alpha}\|^2=\mathcal{O}(n^{-\frac{4\mu p}{4\mu p
+ 2p +1}}).$$ \vskip .1in

Analogous results are obtained in the case of Hilbert scales
(\cite{cohen},\cite{Jmi}).

%Interpreting this rate in the statistical literature reads $s =
%2\mu p$, where the optimal rates are of order
%$\mathcal{O}(n^{-\frac{2s}{2s + 2p +1}})$ , see for example
%\cite{Cavalier}.
\vskip .1in

Adaptive model selection is a technique which penalizes the
regularization parameter, in such a way that we choose
$\hat{f}_{m_0, \alpha_{\hat{k}}}$ by minimizing \vskip .1in

$$ {\rm arg}\min_{k \in \mathcal{K},\, f \in F_{m}} \left(\,\|R_{\alpha_k}(y_m-A_mf)\|^2
+ {\rm pen}(\alpha_k)\, \right)$$
where $$\hat{k}={\rm arg}\min_{k \in
\mathcal{K}}\left(\,\|R_{\alpha_k}(y_m-A_mf)\|^2  + {\rm
pen}(\alpha_k)\, \right)$$ and

$${\rm pen}(\alpha_k)= r\sigma^2(1+L_k)
[Tr(R_{\alpha_k}^t R_{\alpha_k})+\rho^2( R_{\alpha_k})],$$ with
$r>2$ and $L_k$ is a sequence which is incorporated in order to
control the complexity of the set
$\mathcal{K}=\{1,2,\ldots,k_n\}$, of all possible index up to
$k_n$. Here $\rho^2(B) = \rho (B^t B)$ is the spectral radius of
the selfadjoint  operator $B^tB$ for any square matriz $B$, which
is defined by

$$\rho^2(B)= \frac {1}{n} \max_{j\in m} b_j$$
for $b_j$ eigenvalues of $B^tB$. \vskip .1in\vskip .1in

Thus, $\hat{k}$ is selected by minimizing

\be{kestimado}{\rm arg}\min_{k \in
\mathcal{K}}\left(\,\|R_{\alpha_k}(y_m-A_mf)\|^2  +
\frac{r\sigma^2(1+L_k)}{n} \bigg[\sum_{j\in
m}Q^2(\lambda_j)\lambda_j + \max_{j\in
m}Q^2(\lambda_j)\lambda_j\bigg]\, \right).\ee\vskip .1in

The strategy as proposed in this article automatically provides
the optimal order of accuracy. The regularized estimator has a
rate of convergence less or equal than the best rate achieved by
the best estimator for a selected model. We have the following
result,\vskip .2in

\begin{thm}\label{teogenial} For any $f\in F_m$ and any  $\alpha_k$ the following inequality holds true for
$d$ a positive constant that depends on r (as in Lemma
\ref{lema1}),

\be{fin} \BBe \|\tilde{f}_m - \hat{f}_{\alpha_{\hat{k}}}\|^2  \le
\frac{1}{(1-\nu)} \inf_{k \in \mathcal{K}}[C (1+\nu) \|\tilde{f}_m
- f_{\alpha_k}\|^2 + 2{\rm pen}(\alpha_k)] + \frac{C_1(d)}{n} \ee
where $C_1(d)= 4 \, \sigma^2  \sum_{k} \frac{n
\rho^2(R_{\alpha_k})}{d} \left[\sqrt{d\, r \, L_k \left[\frac{
Tr(R_{\alpha_k}^tR_{\alpha_k})}{\rho^2(R_{\alpha_k})}+ 1\right]}
+1 \right] \ e^{ -\sqrt{d \,r \,L_k \left[\frac{
Tr(R_{\alpha_k}^tR_{\alpha_k})}{\rho^2(R_{\alpha_k})}+1 \right]}}.
$
\end{thm}
\vskip .1in

\begin{rem}
An important issue is that equation \eqref{eq:fin} is non
asymptotic. The goodness of fit of the estimator is defined by
trace, $Tr(R_{\alpha}^tR_{\alpha})$, and spectral radius,
$\rho^2(R_{\alpha})$. Also, the estimator is optimal in the sense
that the adaptive estimator achieves the best rate of convergence
among all the regularized estimators.
\end{rem}
\vskip .1in

\begin{rem}
Remark that under our assumptions, namely that the basis is
orthonormal for the fixed design, both $n\rho^2(R_k)$ and
$Tr(R_k^tR_k)/\rho^2(R_k)$ do not depend on n.
\end{rem}
\vskip .1in

\begin{proof}

For any $f_{\alpha_k}$ and any $k \in \BBn$

\bean &\|R_{\alpha_{\hat{k}}}(y_m -
A_m\hat{f}_{\alpha_{\hat{k}}})\|^2 + pen(\alpha_{\hat{k}}) \leq
\|R_{\alpha_k}(y_m - A_m f_{\alpha_k})\|^2 + pen(\alpha_k) \eean
and \bean & \|R_{\alpha_k}(y_m - A_m f_{\alpha_k})\|^2=
\|R_{\alpha_k}A_m(\tilde{f} - f_{\alpha_k})\|^2 + 2 \langle
R_{\alpha_k} A_m (\tilde{f} - f_{\alpha_k}),
R_{\alpha_k}\Pi_{Y_m}^n\varepsilon \rangle + \|R_{\alpha_k}
\Pi_{Y_m}^n\varepsilon\|^2 \eean \vskip .1in

Thus, following standard arguments we have

\bean &&\|R_{\alpha_{\hat{k}}}A_m(\tilde{f}-\hat{f}_{\alpha_{\hat{k}}})\|^2\\
&\leq& \|R_{\alpha_k}A_m(\tilde{f}-f_{\alpha_k})\|^2 - 2<
R_{\alpha_{\hat{k}}}
A_m(\tilde{f}-\hat{f}_{\alpha_{\hat{k}}}),R_{\alpha_{\hat{k}}}\Pi_{Y_m}^n\varepsilon>\\
&&+2<R_{\alpha_{k}}
A_m(\tilde{f}-f_{\alpha_{k}}),R_{\alpha_k}\Pi_{Y_m}^n\varepsilon>-
\|R_{\alpha_{\hat{k}}}\Pi_{Y_m}^n\varepsilon\|^2+
\|R_{\alpha_k}\Pi_{Y_m}^n\varepsilon\|^2 +{\rm pen}(\alpha_k)+{\rm
pen}(\alpha_{\hat{k}}). \eean \vskip .1in

Let $0<\nu<1$. Since the algebraic inequality $2ab\leq \nu
a^2+\frac{1}{\nu}b^2$ holds for all $a,b \in \BBr$, we find that

\bean &&
(1-\nu) \|R_{\alpha_{\hat{k}}}A_m(\tilde{f}-\hat{f}_{\alpha_{\hat{k}}})\|^2\\
&\leq& (1+\nu)\|R_{\alpha_k}A_m(\tilde{f}-f_{\alpha_k})\|^2 + 2
{\rm pen}(\alpha_k) + 2
\sup_{\alpha_k}\{\frac{1}{\nu}\|R_{\alpha_k}\Pi_{Y_m}^n\varepsilon\|^2
- {\rm pen}(\alpha_k)\}, \eean holds for any $k$ and $f_{\alpha_k}
\in F_m$.\vskip .1in

On the other hand, using that is $1\leq \|R_{\alpha_k}A\| \leq C$,
we have that for any $f_{\alpha_k} \in F_{m_0}$ and any $k \in
\BBn$,

\bean && (1-\nu) \|\tilde{f}_m - \hat{f}_{\alpha_{\hat{k}}}\|^2 \le \,C (1+\nu) \|\tilde{f}_m - f_{\alpha_{k}}\|^2\\
&+& 2\,{\rm pen}(\alpha_k) + 2 \, C_1 \,
\sup_{\alpha_k}\,\{\|R_{\alpha_k}\Pi_{Y_m}^n\varepsilon \|^2 -
{\rm pen}(\alpha_k)\}.\eean

\vskip .1in

The proof then follows directly from the following technical lemma
(\cite{Bouquet},\cite{Jmi}) which characterizes the supremum of an
empirical process by the regularization family.\vskip .2in

\begin{lem} \label{lema1}
Let $\eta(A)=\sqrt{\varepsilon^t A^t A \varepsilon} = \|A
\varepsilon\|.$ Then,
%\begin{itemize}
%\item
there exists a positive constant $d$ that depends on $r/2$
such that the  following inequality holds \bea{desprob}
&&P(\eta^2(A)\ge \sigma^2[Tr(A^tA)+ \rho(A^tA)]r/2(1+L)+\sigma^2
u)\\ \nonumber &&\le \exp\{-\sqrt{d (1/\rho(A^tA) u+ r/2 L
[Tr(A^tA)/\rho(A^t A)+1] )}\}.\eea %\vskip .2in \item Set $k_1=
%d/(\rho(A^tA))$ and $k_2=d r/2 L [Tr(A^tA)/\rho(A^tA)+1]$. There
%exists a constant $C_q$, which depends only on $q$, such that,
%\bea{desesp}&& \BBe
%[\eta^2(A)-\sigma^2[Tr(A^tA)+ \rho(A^tA)]r/2(1+L) ]_+^q\\
%\nonumber &&\le C_q \sigma^{2q} k_1^{-q} [ k_2^{q-1/2}+k_2^{q-1}]
%e^{-\sqrt{  k_2}} \eea  holds.
%\end{itemize}
\end{lem}
\vskip .2in

With the above notation,
$$\eta(R_{\alpha_k})=\|R_k \varepsilon_m\|$$
where $\varepsilon_m=\Pi_{Y_m}^n\varepsilon$. \vskip .1in Now,
with this lemma we have

\bean && P(\sup_{\alpha_k}\|R_{\alpha_k}\varepsilon_m\|^2-{\rm
pen}(\alpha_k)
>\sigma^2 x)\\
& \leq & \sum_k P[\eta^2(R_{\alpha_k})\geq r\sigma^2(1+L_k)
[Tr(R_{\alpha_k}^tR_{\alpha_k}) + \rho^2(R_{\alpha_k})] + \sigma^2
x]\\ & \leq & \sum_k \exp\{-\sqrt{d (1/\rho^2(R_{\alpha_k})x +
rL_k [Tr(R_{\alpha_k}^t R_{\alpha_k})/\rho^2(R_{\alpha_k})+1] )}\}
\eean

Since for $X$ positive $\BBe X = \int_0^{\infty} \, P(X > u)\,
du$, we then have that

\bean && \BBe [\sup_{\alpha_k}\|R_{\alpha_k}\varepsilon_m\|^2-{\rm
pen}(\alpha_k)] = \int_0^{\infty} P[\sup_{\alpha_k}
\|R_{\alpha_k}\varepsilon_m\|^2 -{\rm pen}(\alpha_k)\geq x] \,dx\\
& = &
\sigma^2\int_0^{\infty}P[\sup_{\alpha_k}\|R_{\alpha_k}\varepsilon_m
\|^2 -{\rm pen}(\alpha_k) \geq
\sigma^2 u] \,du\\
& \leq & \sigma^2 \sum_k \int_0^{\infty}\exp \{-\sqrt{k_1 u + k_2
}\} \,du.  \eean where $k_1=d/\rho^2(R_{\alpha_k})$ and $k_2= d r
L_k [Tr(R_{\alpha_k}^T R_{\alpha_k})/\rho^2(R_{\alpha_k})+1]$.
\vskip .1in

Let $w=k_1 u + k_2,$ then
\begin{align*}
\BBe [\sup_{\alpha_k} \|R_{\alpha_k}\varepsilon_m \|^2 -{\rm
pen}(\alpha_k)] & \leq \sigma^2 \sum_k
\int_{k_2}^{\infty} \frac{1}{k_1} \exp \{-\sqrt{w}\} \,dw\\
& = \sigma^2 \sum_k \frac{2}{k_1}[-\sqrt{k_2} + 1] \exp
\{-\sqrt{k_2}\}
\end{align*}

Finally, we have the desired result.
\begin{align*}
\BBe \|\tilde{f}_m - \hat{f}_{\alpha_{\hat{k}}}\|^2  \le
\frac{1}{(1-\nu)} \inf_{k \in \mathcal{K}}[C (1+\nu) \|\tilde{f}_m
- f_{\alpha_k}\|^2 + 2{\rm pen}(\alpha_k)] + \frac{C_1(d)}{n}
\end{align*}
where $C_1(d)= 4 \, \sigma^2  \sum_{k} \frac{n
\rho^2(R_{\alpha_k})}{d} \left[\sqrt{d\, r \, L_k \left[\frac{
Tr(R_{\alpha_k}^tR_{\alpha_k})}{\rho^2(R_{\alpha_k})}+ 1\right]}
+1 \right] \ e^{ -\sqrt{d \,r \,L_k \left[\frac{
Tr(R_{\alpha_k}^tR_{\alpha_k})}{\rho^2(R_{\alpha_k})}+1 \right]}},
$
\end{proof}
%%% ----------------------------------------------------------------------
\vskip .3in

\section{Regularization by iterative methods }\label{RIM} \vskip .2in

Iterative regularization methods, are very competitive methods for
linear inverse problems. In iterative regularization, one picks an
initial guess $f_0$ for the unknown $\tilde{f}$, and then one
iteratively constructs updated approximations via a regularization
scheme. The regularization parameter associated with iterative
regularization is thus the $\lq \lq$stopping point\rq \rq of the
iterative sequence, and an important part of the mathematical
theory is the development of stopping criteria for terminating the
iteration. In other words, the iteration index plays the role of
the regularization parameter $\alpha$, and the stopping criteria
plays of the parameter selection method. \vskip .3in

%%% ----------------------------------------------------------------------

\subsection{Descent Methods for Linear Inverse Problems}$ $ \vskip .2in

As an example of iterative regularization, we consider descent
methods. Descendent methods have become quite popular in the last
years for the solution of linear inverse problems and for
nonlinear inverse problems \cite{Gilyazov}. In this subsection we
consider two examples.
%These results say how the dimension of the
%subspace $d_m$ must be chosen, in terms of the regularization
%parameter $\alpha$. They also say, that under additional technical
%conditions these results are good, in the sense they achieve
%optimal rates. The ill posedness, as measured by the sequence
%${\lambda_j}$ must be considered in this control.
\vskip .1in

As an approximation of $\tilde{f}_m$ we will choose $f_{m,\alpha}$
such that

\be{model1} f_{m,\alpha} = [I-A_m^*A_m Q_{\alpha}(A_m^*A_m)]f_{0}
+ Q_{\alpha}(A_m^*A_m)A_m^*\Pi_{Y_m}^n y\ee where $f_{0} \in F_m$
is an initial approach and this $f_{0} \in
\mathcal{N}(A_m)^{\bot}$ \cite{Gilyazov}. \vskip .1in

Most iterative methods for approximating $\tilde{f}$ are based on
a transformation of the normal equation into equivalent fixed
point equations like

$$f = f + A_m^* (A_m f -y) $$

If $\|A_m\|^2 < 2$ then the corresponding fixed point operator $I
- A_m^* A_m$ is nonexpansive and one may apply the method of
successive approximations. It must be emphasized that $I-A_m^*A_m$
is no contradiction if our inverse problem is ill-posed, since the
spectrum of $A_m^*A_m$ clusters at the origin.\vskip .3in
%%% ----------------------------------------------------------------------

\subsection{Landweber iteration}$ $ \vskip .2in

In this subsection we presented the well-known Landweber
iteration, which arises from converting the necessary conditions
for minimizing \eqref{eq:2} into a fixed point iteration. Much
development in the last few years has taken place in advancing the
theory of Landweber iteration for linear and nonlinear inverse
problems. \vskip .1in

Using the terminology of the last sections, we introduce the
function

\be{ec1} Q_k(\lambda) =  \sum_{j=0}^{k-1}\, (1 - \lambda)^j =
\lambda^{-1}(1-(1-\lambda)^{k})\ee \vskip .1in

We call $Q_k$ the \emph{iteration polynomial} of degree $k-1$.
Associated with it is the polynomial

$$r_k(\lambda) = 1 - \lambda Q_k(\lambda) = (1-\lambda)^k$$
of degree $k$, which is called the \emph{residual polynomial}
since it determines the residual $y - A_m f_{m,k}$.\vskip .1in

Thus, inserting the equation \eqref{eq:ec1} in \eqref{eq:model1}
we obtain recursively,

\be{ecLandweber1} f_{m,k+1} = f_{m,k} -  A_m^*(A_m f_{m,k} -
y_m),\, \,k=0,1,\ldots \ee starting from an initial guess $f_{0}$.
This is a steepest descent method called  the \emph{linear
version} of Landweber's iteration. Each step of the iterative
process \eqref{eq:ecLandweber1} is carried out along the direction
opposite to the direction of the gradient of the quadratic
functional $J(f)$ in \eqref{eq:2}. It is known that there is the
greatest decrease of the functional along this direction.
%Such iterative methods are usually called \emph{gradient methods.} Note
%that, we make the identification $\alpha = k$ for iterative
%methods. This is a first candidate for solving \eqref{eq:2} in an
%iterative way. For nonlinear problems the investigation of an
%iterative method is much more complicated and has mostly to be
%done for each class of methods individually. For more discussion
%on this method and numerous other iterative regularization
%schemes, see \cite{Engl}, \cite{Gilyazov}.\vskip .1in \vskip .1in

%\begin{rem}
%The Landweber iteration, in case of ambiguity, can be supplied
%with and initial guess $f_0 = f^{\star}$ which plays the same role
%as in Tikhonov regularization. Also, if the problem is linear,
%good solutions are obtained for $f_0 = 0$.
%\end{rem}
\vskip .1in

If $\|A_m \| \leq 1 $, we considerer $ \lambda \in (0,1]$ such
that in this interval $\lambda Q_k(\lambda)$ is uniformly bounded
and since $Q_k(\lambda)$ converge to $1/\lambda$ as $k \to \infty$
then according to Theorem  \ref{teo1} the sequences $f_{m,k}$
converge to $\tilde{f}_m$ when $y \in \mathcal{D}(A_m^{\dag})$.
%On other hand, the $k$th iterate depends of space dimension. In the
%following we derive a simple estimate for the error propagation in
%the Landweber iteration.
%Throughout the analysis it will turn out convenient to have
%$\|A_m\| \leq 1$.
If $\|A_m\|$ is not bounded by one, %this were no the case,
then we introduce a relaxation parameter $0 < \tau < \|A_m\|^{-2}$
in front of $A_m^*$ in \eqref{eq:ecLandweber1}, i.e, we would
iterate

\be{ecLandweber2} f_{m,k+1} = f_{m,k} -  \tau A_m^*(A_m f_{m,k} -
y), k=0,1,\ldots \ee

If $\tau \equiv \tau_k$, one can obtain various variants of the
method of steepest descent depending on a choice of the sequence
$\tau_k$. The Landweber iteration \eqref{eq:ecLandweber2} is
usually called \emph{a method of simple iteration}. \vskip .1in

In the following we derive a simple estimate for the error
propagation in the Landweber iteration. We then have the following
result,\vskip .2in

\begin{cor} Let $\tau = 1/ (2 \|A_m\|^{2}) < 1/\lambda_1$. If $y \in \mathcal{R}(A_m)$, then the Landweber iteration
is an order optimal regularization method, i.e,
$$\|\tilde{f}_m- f_{k(m)}\|^2 \leq 2 c_1 k^{-2\mu} + 2 c_2
\frac{\sigma^2}{n} (\tau k)^{(2p+1)/2p},$$ where
$c_1=\rho^2(\frac{\mu}{\tau e})^{\mu}$ and
$c_2=\frac{1}{2p+1}(\frac{2p+1}{2p-1})^{(2p+1)/4p}.$
\end{cor}
\vskip .1in

\begin{proof} To apply Theorem \ref{teo2} we have to study the terms of the bias,
$\BBe\|\tilde{f}_m-R_{\alpha}A_m\tilde{f }\|^2$, and variance
$\BBe\| R_{\alpha}A_m \tilde{f} - f_{k(m)}\|$. By \eqref{eq:ec1}
we have

$$\tilde{f}_m- R_{\alpha}A_m\tilde{f} = (I - A_m^* A_m Q_k(A_m^*A_m))\tilde{f}_m = (I - A_m^*A_m)^k \tilde{f}_m$$
\vskip .2in We have to study the residual polynomial
$r_k(\lambda)=(1-\lambda)^k$ of the Landweber iteration.\vskip
.1in

For $0 \leq \lambda \leq \|A_m\|^2$ the function
$$\lambda^{\mu}|1-\lambda Q_k(\lambda)|$$ assumes its maximum for
$\lambda = \tau^{-1} \mu (\mu + k)^{-1}$.\vskip .1in

Thus, we have
\begin{align*}
\lambda^{\mu}|1-\lambda Q_k(\lambda)| & \leq \max_{0\leq
\lambda\leq \|A_m\|^2} \lambda^{\mu}|1-\lambda Q_k(\lambda)|\\&  <
\frac{\mu^{\mu}}{\tau^{\mu}(\mu + k)^{\mu}} \frac{k^k}{(\mu +
k)^k} \\&<  \bigg(\frac{\mu}{\tau e}\bigg)^{\mu} k^{-\mu}
\end{align*}\vskip .1in

This leads to numbers $\omega_{\mu}(k)$ as introduced in  Theorem
\ref{teo2}

$$ \omega_{\mu}(k) = \bigg(\frac{\mu}{\tau e}\bigg)^{\mu} k^{-\mu}$$

Thus, the term corresponding to the bias is bounded by
$$\|\tilde{f}_m-R_{\alpha}A_m\tilde{f }\|^2 \leq
\rho^2(\frac{\mu}{\tau e})^{2\mu}k^{-2\mu}.$$% \vskip .2in

Next, we establish bounds for the variance term. By assumption,
the singular values satisfy $\lambda_j \approx j^{-2p}.$ Note that
for small values of $\lambda_j$  we have
\begin{center}
 $Q^{2}_k(\lambda) \leq
(\tau k)^{2}$ \, \,$ \forall$ $j>m^{\prime}$
\end{center}
and for big values of $\lambda_j$ ($\lambda_j \approx \lambda_1$)

\begin{center}
$Q^{2}_k(\lambda_j)\leq \lambda_j^{-2}$ \, \,$ \forall$
$j<m^{\prime}$.
\end{center}

\vskip .1in

Consequently,
\begin{align*}
n Tr(Q_k^{2}(A_mA_m^{*})A_mA_m^{*})=
\sum_{j=1}^{m}Q_k^{2}(\lambda_j)\lambda_j & \leq
\sum_{j=1}^{m^{\prime}} \lambda_j^{-1} +
\sum_{j> m^{\prime}}(\tau k)^{2}\lambda_j\\
& \leq \int_0^{m^{\prime}}s^{2p}ds + (\tau k)^{2}
\int_0^{m^{\prime}}s^{-2p}ds
\end{align*}

This suggest searching $m^{\prime} \approx c (\tau k)^{1/2p}$ for
$p>1/2$, where $c=(\frac{2p+1}{2p+1})^{1/4p}$. Hence we have,% \vskip .1in

\begin{align*}
\BBe\| R_{\alpha}A_m \tilde{f} - f_{k(m)}\|^2 & = Tr(Q_k^{2}(A_m
A_m^{*})A_m A_m^{*}) \\ & \leq \frac{c^{2p+1}}{2p+1}\frac{(\tau
k)^{(2p+1)/2p}}{n}.
\end{align*}
%\vskip .2in

Finally, this implies
$$\BBe \|\tilde{f}_m-f_{k(m)}\|^2 \leq 2\, c_1 k^{-2\mu} + 2 \,c_2
\frac{\sigma^2}{n} (\tau k)^{(2p+1)/2p},$$ where
$c_1=\rho^2(\frac{\mu}{\tau e})^{\mu}$ and
$c_2=\frac{1}{2p+1}(\frac{2p+1}{2p-1})^{(2p+1)/4p}$.\vskip .1in
\vskip .1in

\end{proof}

\begin{rem}
Note that under the above inequality is satisfied if the dimension
of the set is such that $d_{m_0} \approx n^{\frac{1}{4\mu p + 2p +
1}}$. Here, the optimal choice of regularization sequence,
depending on $p$ and $\mu$. The optimal rates are of order $\BBe
\|\tilde{f}-f_{k(m)}\|^2 = \mathcal{O}(n^{- \frac{4 \mu p}{4 \mu p
+ 2p +1}})$. Analogous results are obtained in the ill-posed
problem literature, see for example \cite{Cavalier}, where
typically in a Hilbert scale setting optimal rates are of order
$\mathcal{O}(n^{- \frac{ 2s}{2s
 + 2p +1}})$, with $s=2\mu p$.
\end{rem}
\vskip .1in

We are ready to state our main result for the Landweber iteration,
which bounds the mean squared error of the select estimate
$\hat{f}_{\hat{k}}$ basically by the smallest mean squared error
among the estimates $f_k$ plus a remainder term of order $1/n$.
The result follows  from Theorem \ref{teogenial}. \vskip .1in

\begin{cor}\label{corolario2}
Let $\tau = 1/ (2 \|A_m\|^{2}) < 1/\lambda_1$. Next assume
$\hat{k}$ as in \eqref{eq:kestimado} and $d_{m_0}$ as in
\eqref{eq:condition}. If $y \in \mathcal{R}(A_m)$ then for any $f
\in F_m$ and any $k$, the following inequality holds true

\begin{align*}
\BBe \|\tilde{f}_m - \hat{f}_{\hat{k}}\|^2  & \leq
\frac{1}{(1-\nu)} \inf_{k \in \mathcal{K}}\bigg[C (1+\nu)
\|\tilde{f}_m - f_{k}\|^2 + \frac{2r\sigma^2(1+L_k)(c(\tau
k)^{\frac{2p+1}{2p}} + \tau k)}{n}\bigg] \\ & +  \frac{4
\sigma^2}{n} \sum_{k} \frac{\tau k}{d} \left[\sqrt{d\, r \, L_k
\left[c(\tau k)^{1/2p}+ 1\right]} +1 \right] \ e^{ -\sqrt{d \,r
\,L_k \left[c(\tau k)^{1/2p}+1 \right]}},
\end{align*}
for some $C>0$ and
$c=\frac{1}{2p+1}(\frac{2p+1}{2p-1})^{(2p+1)/4p}$.

\end{cor}
\vskip .1in

%\begin{rem}
%The goodness of fit of the estimator is defined by trace,
%$Tr(R_{\alpha}^tR_{\alpha})$, and spectral radius,
%$\rho^2(R_{\alpha})$. The condition $d_{m_0}\geq k^{1/2p}$ follows
%of
%$$\frac{Tr(R_{\alpha}^tR_{\alpha})}{\rho^2(R_{\alpha})}=\mathcal{O}(k^{1/2p})$$
%\end{rem}
%\vskip .1in

\begin{proof}
For fixed $\lambda_j$ and $m$ we have that the terms of the trace
and spectral radius are bounded by the follows expression

\be{traza} Tr(R_k^tR_k)= \sum_{j=1}^{m}Q_k^{2}(\lambda_j)\lambda_j
 \leq \sum_{j=1}^{m^{\prime}} \lambda_j^{-1} + \sum_{j>m^{\prime}}(\tau
k)^{2}\lambda_j\ee

and

\be{radio} \rho^2(R_k^tR_k)= \max_{j \in m} \,
Q_k^{2}(\lambda_j)\lambda_j  \leq \max_{j \leq
m^{\prime}}\lambda_j^{-1} + \max_{j > m^{\prime}}(\tau
k)^{2}\lambda_j .\ee \vskip .1in

Balancing both terms in \eqref{eq:traza} and \eqref{eq:radio}
gives the optimal choice of the trace and the spectral radius,
respectively. Thus, we have

$$Tr(R_k^tR_k)\approx
\frac{1}{2p+1}\bigg(\frac{2p+1}{2p-1}\bigg)^{\frac{(2p+1)}{4p}}
\frac{(\tau k)^{\frac{2p+1}{2p}}}{n}$$
 and
$$\rho^2(R_k^tR_k) \approx \frac{\tau k}{n}$$

Note that the penalization term is roughly proportional to

$$\frac{1}{n} \bigg[ \frac{1}{2p+1}\bigg(\frac{2p+1}{2p-1}\bigg)^
{\frac{(2p+1)}{4p}}(\tau k)^{\frac{2p+1}{2p}} + \tau k\bigg ] $$

On the other hand

$$\frac{Tr(R_{\alpha}^tR_{\alpha})}{\rho^2(R_{\alpha})}=
\frac{1}{2p+1}\bigg(\frac{2p+1}{2p-1}\bigg)^{(2p+1)/4p} (\tau
k)^{1/2p}.$$ \vskip .1in

The result then follows directly from Theorem \ref{teogenial}.
\end{proof}

\vskip .3in
%%% ----------------------------------------------------------------------

\subsection{Nonlinear multistep iterative process}$ $ \vskip .2in

Many approximate methods widely used in practice are nonlinear. We
cite a important example of nonlinear approximate method. We
considerer a nonlinear multistep iterative  process, which have
error residual
$$1-\lambda Q_k(\lambda)= \prod_{i=1}^k (1-\tau_{ik}^{-1} \lambda)$$
with $\tau_{ik}=\tau_{ik}(f_0,A,y)>0,\, 0<\tau_{1k}\leq
\tau_{2k}\ldots \leq \tau_{kk}\leq \lambda_{1}$. Then for $\lambda
> 0$, $Q_k(\lambda)$ have the following representation
$$Q_k(\lambda)=\lambda^{-1}[1-\prod_{i=1}^k
(1-\tau_{ik}^{-1} \lambda)]$$ \vskip .1in

The following corollary is established.\vskip .2in

\begin{cor}\label{corolario3} Let $\tau_{ik}=\tau_{ik}(f_0,A,y)>0,$\, with $0<\tau_{1k}\leq
\tau_{2k}\ldots \leq \tau_{kk}\leq \lambda_{1}$. If $y \in
\mathcal{R}(A_m)$, then the nonlinear multistep iterative process
is an order optimal regularization method, i.e,
$$\BBe \|\tilde{f}_m- f_{k(m)}\|^2 \leq 2\, c_1\, (\sum_{i=1}^{k}\tau^{-1}_{i_k})^{-2\mu} + 2\,
c_2\, \frac{\sigma^2}{n}\,
(\sum_{i=1}^{k}\tau^{-1}_{i_k})^{(2p+1)/2p},$$ where
$c_1=\rho^2\mu^{\mu}(\mu+1)^{-1}$ and
$c_2=\frac{1}{2p+1}(\frac{2p+1}{2p-1})^{(2p+1)/4p}.$

\end{cor}
\vskip .1in

\begin{proof}

As before, we investigate the behavior of the bias and the
variance.\vskip .1in

In the relation
$$\lambda^{\mu}(1-\lambda Q_k(\lambda))= \lambda^{\mu} \prod_{i=1}^k (1-\tau_{ik}^{-1} \lambda)$$
the least upper can not be reached at the points $\lambda=0$ and
$\lambda=\tau_{1k}$, since the estimated function is not equal to
zero identically. \vskip .1in

On the other hand

\be{999}[1-\lambda Q_k(\lambda)]'[\lambda Q_k(\lambda)-1]^{-1} =
\sum_{i=1}^k \frac{1}{\tau_{ik}-\lambda}\ee Since the function in
the right-hand  of \eqref{eq:999} does not decrease  as a function
of $\lambda$ on the half-interval $[0,\tau_{1k})$ then, the
estimates of the Theorem \ref{teo3} are valid.\vskip .1in

Thus, for $0\leq \lambda\leq \tau_{1k}$, we have

    \begin{align*}
    \sup_{0\leq \lambda\leq \tau_{1k}}|Q_k(\lambda)|  & = Q_k(0) =  \sum_{i=1}^{k}\tau^{-1}_{i_k}
    \end{align*}
and
    \begin{align*}
    \sup_{0\leq \lambda\leq \tau_{1k}}\lambda^{\mu}|1-\lambda Q_k(\lambda)| & < \mu^{\mu}(\mu+1)^{-1} (\sum_{i=1}^{k}\tau^{-1}_{i_k})^{-\mu}
    \end{align*}
\vskip .1in

Note that
$\omega_{\mu}(k)=(\sum_{i=1}^{k}\tau^{-1}_{i_k})^{-\mu}.$ Thus,
the bias  is bounded by

$$\|\tilde{f}_m-R_{\alpha}A_m\tilde{f }\|^2 \leq c_1
(\sum_{i=1}^{k}\tau^{-1}_{i_k})^{-2 \mu},$$ where
$c_1=\rho^2\mu^{2\mu}(\mu+1)^{-2}$. \vskip .1in \vskip .1in

On the order hand, it is not difficult to see that
\begin{align*}
Tr(Q_k^{2}(A_mA_m^{*})A_mA_m^{*}) & \leq \frac{1}{n}\bigg[\sum_{1
\leq j \leq m^{\prime}}j^{2p} +
(\sum_{i=1}^{k}\tau_{ik}^{-1})^{2}\sum_{j>
m^{\prime}}j^{-2p}\bigg] \leq
\frac{1}{n}\bigg[\frac{{m^{\prime}}^{2p+1}}{2p+1} +
(\sum_{i=1}^{k}\tau_{ik}^{-1})^{2}\frac{{m^{\prime}}^{-2p+1}}{2p-1}\bigg]
\end{align*}
\vskip .1in

This suggest searching

$$m^{\prime} \approx c (\sum_{i=1}^{k}\tau^{-1}_{i_k})^{1/2p}$$ with
$c=(\frac{2p+1}{2p+1})^{1/4p}$.\vskip .2in

Thus, we have that the term variance is bounded by
$$\BBe\| R_{\alpha}A_m \tilde{f} - f_{k(m)}\| = \sigma^2 Tr(Q_k^2(A_m^*A_m)A_m^*A_m) \leq c_2 \frac{\sigma^2}{n} (\sum_{i=1}^{k}\tau^{-1}_{i_k})^{(2p+1)/2p}$$
where $c_2=\frac{1}{2p+1}(\frac{2p+1}{2p-1})^{(2p+1)/4p}$

\vskip .2in

Finally we have
$$\BBe \|\tilde{f}_m-f_{k(m)}\|^2 \leq 2\, c_1\, (\sum_{i=1}^{k}\tau^{-1}_{i_k})^{-2\mu} + 2\,
c_2\, \frac{\sigma^2}{n}\,
(\sum_{i=1}^{k}\tau^{-1}_{i_k})^{(2p+1)/2p}$$

Balancing the bias and variance terms gives the optimal choice

 $$\BBe \|\tilde{f}_m-f_{k(m)}\|^2 = \mathcal{O}(n^{- \frac{4 \mu p}{4 \mu p + 2p
 +1}}).$$
\end{proof}

\vskip .1in

We have the following result.\vskip .1in
\begin{cor}\label{corolario4}
Let $\tau_{ik}$ be as in corollary \ref{corolario3} . Next assume
$\hat{k}$ as in \eqref{eq:kestimado} and $d_{m_0}$ as in
\eqref{eq:condition}. If $y \in \mathcal{R}(A_m)$ then for any $f
\in F_m$ and any $k$, the following inequality holds true

\begin{align*}
\BBe \|\tilde{f}_m - \hat{f}_{\hat{k}}\|^2  & \leq
\frac{1}{(1-\nu)} \inf_{k \in \mathcal{K}}\bigg[C (1+\nu)
\|\tilde{f}_m - f_{k}\|^2 +
\frac{2r\sigma^2(1+L_k)(c(\sum_{i=1}^{k}\tau^{-1}_{i_k})^{\frac{2p+1}{2p}}
+ \sum_{i=1}^{k}\tau^{-1}_{i_k})}{n}\bigg] \\ & +  \frac{4
\sigma^2}{n} \sum_{k} \frac{\sum_{i=1}^{k}\tau^{-1}_{i_k}}{d}
\left[\sqrt{d\, r \, L_k
\left[c(\sum_{i=1}^{k}\tau^{-1}_{i_k})^{1/2p}+ 1\right]} +1
\right] \ e^{ -\sqrt{d \,r \,L_k \left[c
(\sum_{i=1}^{k}\tau^{-1}_{i_k})^{1/2p}+1 \right]}},
\end{align*}
for some $C>0$ and
$c=\frac{1}{2p+1}(\frac{2p+1}{2p-1})^{(2p+1)/4p}$.
\end{cor}
\vskip .1in
\begin{proof}
First observe that

$$Tr(R_k^tR_k)\approx
\frac{1}{2p+1}\bigg(\frac{2p+1}{2p-1}\bigg)^{\frac{(2p+1)}{4p}}
\frac{(\sum_{i=1}^{k}\tau^{-1}_{i_k})^{\frac{2p+1}{2p}}}{n}$$
 and
$$\rho^2(R_k^tR_k) \approx \frac{
\sum_{i=1}^{k}\tau^{-1}_{i_k}}{n}$$

Consequently

$$\frac{Tr(R_k^tR_k)}{\rho(R_k)} \approx
\frac{1}{2p+1}\bigg(\frac{2p+1}{2p-1}\bigg)^{\frac{(2p+1)}{4p}}
(\sum_{i=1}^{k}\tau^{-1}_{i_k})^{\frac{1}{2p}} $$

\vskip .1in

Note that both $n\rho^2(R_k)$ and $Tr(R_k^tR_k)/\rho^2(R_k)$ do
not depend on n. The proof then follows directly from theorem
\ref{teogenial}.

\end{proof}

\vskip .3in
%%% ----------------------------------------------------------------------

%\section*{Acknowledgments}

%%% BIBLIOGRAFIA

%%% FIN DE LA BIBLIOGRAFIA
%%% ----------------------------------------------------------------------

\end{document}